\documentclass[a4paper,12pt]{article}
\usepackage{stmaryrd}
\usepackage{bbm}
\usepackage{amsfonts}
\usepackage{amsmath}
\usepackage{amssymb}
\usepackage{mathrsfs}
\usepackage{accents}
\usepackage{amscd}

\newtheorem{atheorem}{\bf \temp}[section]
\newtheorem{thm}[atheorem]{Theorem}

\newtheorem{lem}[atheorem]{Lemma}
\newtheorem{prop}[atheorem]{Proposition}
\newtheorem{de}[atheorem]{Definition}
\newtheorem{rem}[atheorem]{Remark}

\numberwithin{equation}{section}

\title{\textbf{Partial regularity of Solutions of  Navier-Stokes equations }}
\author{ Xixia Ma  \footnote{Email: kfmaxixia@163.com} \ \ \ \ \ \ \
 \\}
\date{}

\begin{document}

\maketitle

\textbf{Abstract.} In this paper, we study the singular set of 3-dimensional  Navier-Stokes equations. Under the condition$\frac{1}{R^{\frac{3s}{q}+2-s}}\int^{R^{2}}_{0}(\int_{B_{R}}|u|^{q}dx)^{\frac{s}{q}}ds <C,$ for $(q,s)\in\{(2,5),(5,2)\},$ we use the backward uniqueness of parabolic equations to show that the Hausdorff dimension of the singular set is less than 1.

\begin{center}
\item\section{Introduction}

\end{center}
$$$$

Let $\Omega$ be a domain in $\mathbb{R}^{3}$ with smooth boundary $\partial\Omega.$ On the space-time cylinder $\Omega\times(0,\infty),$ we consider the incompressible Navier-Stokes equations in three dimensional space with unit viscosity,
\begin {equation}
\left\{\begin{array}{l}
\partial_{t}u-\triangle u+u\cdot\nabla u+\nabla p=0,\quad t>0,\\
\texttt{ div} u=0,\\
u|_{\partial\Omega}(x,t)=0,\\
u(0,x)=u_{0}(x).
\end{array}\right.
\end{equation}
The velocity field $u=(u_{1},u_{2},u_{3}):\Omega\times(0,\infty)\rightarrow\mathbb{R}^{3},$ and $p(x,t):\Omega\times(0,\infty)\rightarrow\mathbb{R}$ is the pressure.  It is a long standing open question to determine if solutions with large smooth initial data of finite energy remain regular for all time.

In this paper, we  consider the special class of solutions which are suitable weak solutions. The definition of suitable weak solutions is introduced in [2] as follows.
 \begin{de}let $\Omega $ be a open set in $\mathbb{R}^{3}$. We say that a pair $u$ and $p$ is a suitable weak solution to the Navier-Stokes equations on the set $\Omega\times(-T_{1},T)$ if it satisfies the conditions :

 $ \textbf{i},$ \begin{equation}
   u\in L^{2,\infty}(\Omega\times(-T_{1},T))\cap L^{2}(-1,0;H^{1}(\Omega)) ,  p\in L^{\frac{3}{2}}(\Omega\times(-T_{1},T));
    \end{equation}

  $ \textbf{ii}$, \quad $u$ and $p$ satisfy the Navier-Stokes equations in the  distribution sense;

   $\textbf{iii}$, \quad $u$ and $p$ satisfy the local energy inequality
   \begin{equation}
    \int_{\Omega}\varphi|u(x,t)|^{2}+2\int_{\Omega\times(-T_{1},t)}\varphi|\nabla u|^{2}dxdt'
    \leq$$ $$\int_{\Omega\times(-T_{1},t)}(|u^{2}(\triangle\varphi+\partial_{t}\varphi)+u\cdot\nabla\varphi(|u|^{2}+\nabla2p))dxdt' \end{equation}
 for a.a. $t\in(-T_{1},T)$ and for all nonnegative functions $\varphi\in C^{\infty}_{0}(\mathbb{R}^{3}). $\end{de}

There are lots of important papers that contribute to the regularity problem of
suitable weak solutions to the Navier-Stokes equations and there are many good survey papers
and books. Hence, we only list some of them. Scheffer [8, 9] introduced partial regularity for
the Navier每Stokes system. Caffarelli, Kohn and Nirenberg [2] further strengthened Scheffer＊s
results. Lin [7] gave a new short proof for the result of Caffarelli, Kohn and Nirenberg. Ladyzhenskaya
and Seregin [10] investigated partial regularity. Choe and Lewis [3] studied singular
set by using a generalized Hausdorff measure. Escauriaza, Seregin, and $\breve{S}$ver$\acute{a}$k [1] proved
the critical case of the so-called Ladyzhenskaya-Prodi-Serrin condition based on the backward unique
continuation theory for parabolic equations. Gustafson, Kang, and Tsai [4] generalize several
previously known criteria. Here we state one of the main results of the theory of suitable weak solutions as follows.
\begin {lem}(see [1]) There exist absolute positive constants $\varepsilon_{0}$ and $c_{0k}, k = 1,2,...,$
with the following property. Assume that a pair u and p is a suitable weak
solution to the Navier-Stokes equations in Q and satisfies the condition
$$\int_{Q}|u|^{3}+|p|^{\frac{3}{2}}dxdt<\varepsilon_{0}$$

Then, for any natural number$ k, \nabla^{k-1}u $ is H$\ddot{o}$lder continuous in $\bar{Q}(\frac{1}{2})$ and
the following bound is valid:
$$\max_{z\in Q(\frac{1}{2})}\nabla^{k-1}u\leq c_{0k}.$$
\end {lem}
\begin {thm}(see[2]) For any suitable weak solution of the Navier-Stokessystem on an
open set in space-time, the associated singular set satisfies  $\mathfrak{H}^{1}(S)=0.$
\end {thm}

\begin {rem} Let us mention, from a physical point of view, the result of Theorem 1.3(Caffarelli,Kohn, and Nirenberg[1982]) gave an answer about Jean Leray's conjecture concerning the appearance of singularities in 3-dim turbulent flow, that is , if there exists a singular set which is a fractal set, then the occurrence of smooth line vortices is not possible. Furthermore, this powerful mathematical result leaves room for a tremendously complex set of singularites, and we remain far from closing the issues raised by Leray's conjecture\footnote{Leray's conjecture: turbulence on Navier-Stokes equations is due to the formation of point or "line vortices" on which some component of the velocity becomes infinite. }.
 \end {rem}

To enable dealing with his conjecture, Leray suggested the concept of weak, nonclassical solutions to the Navier-Stokes equations(1.1), and this has become the starting point of the mathematical theory of the Navier-Stokes equations to this day. However, even today, J.Leray's conjecture concerning the appearance of singularities in 3-dimensional turbulence flows has been neither proved nor disproved. In this paper, we try to improve the results of Theorem 1.3 through the following theorem.
\begin {thm}
  For any suitable weak solution $(u,p)$ of the Navier-Stokes system on an open set in space-time , and $(u,p)$ satisfies the following condition:
  \begin {equation}
  \frac{1}{R^{\frac{3s}{q}+2-s}}\int^{R^{2}}_{0}(\int_{B_{R}}|u|^{q}dx)^{\frac{s}{q}}ds <C
  \end {equation}
  where $ C $is a larger absolute constant, and $(q,s)\in\{(2,5),(5,2)\}.$ Meanwhile, $ p$  satisfies \begin {equation}
  \frac{1}{R^{2}}\int^{R^{2}}_{0}(\int_{B_{R}}|p|^{\frac{3}{2}}dx)ds <C.
  \end {equation}
   Then the parabolic Hausdorff dimension of the associated singular set is at most $\frac{10-m}{m+2},$ for any $ m \in (4,5).$
  \end {thm}

  The idea of our proof is from the scaling invariant of Navier-Stokes equations,using the blow-up procedure, we use the backward uniqueness results of parabolic equations. The article is organized as follows. Some auxiliary results are given in section 2,  and we will give the proof of our main theorem in the last section.

  \begin{center}
\item\section{A new $\varepsilon-$ Regularity Criterion }
\end{center}

In this section, first, we give some preliminary. Furthermore, using the backward uniqueness property of parabolic equations, we prove a new $\varepsilon-$ regularity criterion of Naiver-Stokes equations.

\begin {thm} (Backward uniqueness for Heat Operator) We consider a vector-valued function $u: (\mathbb{R}^{n}\setminus B(R))\times[0,T]\rightarrow \mathbb{R}^{n},$ assume $u$ satisfies the following conditions:

(a) \begin{equation}
|\partial_{t}u+\triangle u|\leq c_{1}(|\nabla u|+|u|)\quad \textsf{in}\quad(\mathbb{R}^{n}\setminus B(R))\times[0,T]
\end{equation}
for some $c_{1}>0 ;$

(b) \begin{equation}
u(\cdot,0)=0 \quad \textsf{in}\quad(\mathbb{R}^{n}\setminus B(R));
\end{equation}

 (c) \begin{equation}
 |u(x,t)|\leq \exp^{M|x|^{2}}
 \end{equation}
 for all $(x,t)\in(\mathbb{R}^{n}\setminus B(R))\times[0,T]$ and for some $ M>0;$
 
 (d) $ u$ and distributional derivatives $\partial_{t}u,\nabla^{2}u$ are square integrable over bounded subdomains of $(\mathbb{R}^{n}\setminus B(R))\times[0,T].$
 Then $u\equiv0$ in $(\mathbb{R}^{n}\setminus B(R))\times[0,T].$
 
\end {thm}

The proof of Theorem 2.1 comes from in [],so we omit the proof.
Note the estimates in the following lemma is scaling invariance.

\begin{lem} Assume $(u,p)$ is a Leray-Hopf solution of (1.1), if $ u $ satisfies 
 $$ \frac{1}{R^{\frac{3s}{q}+2-s}}\int^{R^{2}}_{0}(\int_{B_{R}}|u|^{q}dx)^{\frac{s}{q}}ds <C ,$$ where $ C $is an absolute constant, and $(q,s)\in\{(2,5),(5,2)\}.$ Then we have $\partial_{t}u,\nabla p,\nabla^{2} u\in L^{\frac{2m}{m+2}}(Q_{T})$ for any $ m\in(4,5)$.\end{lem}

  \textbf{proof}: First by H$\ddot{o}$lder inequality and Sobolev embedding,
  $$
    \|u\cdot\nabla u\|_{L^{\frac{2m}{m+2}}}\leq\|u\|_{L^{m}}\|\nabla u\|_{L^{2}}$$

   Choose $\phi\in C^{\infty}_{0}(\mathbb{R}^{3})$ and $ div \phi =0$, we have
   \begin {equation}
   \begin{array}{l}
  (\partial_{t}u,\phi)=-(u\nabla u,\phi)-(\nabla u,\nabla\phi)\\
                     \leq\|u\|_{L^{2}}\|\nabla u\|_{L^{2}}\|\phi\|_{L^{\infty}}+\|\nabla u\|_{L^{2}}\|\nabla\phi\|_{L^{2}}\\
                   \leq(\|u\|_{L^{2}}+\|u\|_{L^{2}}\|\nabla u\|_{L^{2}})\|\phi\|_{H^{\frac{3}{2}}}
                    \end{array}
                    \end{equation}.
we know that $$
\|u\|_{L^{m}(Q_{\frac{1}{2}})}\leq\|u\|^{\theta}_{L^{5}((-\frac{1}{4},\frac{1}{4});L^{2}(B_{\frac{1}{2}}))}
\|u\|^{1-\theta}_{L^{2}((-\frac{1}{4},\frac{1}{4});L^{5}(B_{\frac{1}{2}}))}
$$
with $\theta=\frac{2(5-m)}{3m}.$
                   Hence $\partial_{t}v\in L^{\frac{2m}{m+2}}(Q_{T})$,since $ m\in(4,5)$. In the following , we show $\nabla p\in L^{\frac{2m}{m+2}}(Q_{T})$ ,then it is easy to check $\nabla^{2} v\in L^{\frac{2m}{m+2}}(Q_{T}$.

                   In fact , let $ f=\partial_{t}v-\triangle v $ ,then first it is obtained that
                   $f\in L^{2}(0,T;H^{-\frac{3}{2}}_{0})  $ as mentioned above.

                   And then we know
\begin {equation}
 \left\{\begin {array}{l}
   \textsf{ div} f =0, \\
   \ast df= \ast d(v\cdot\nabla v)
   \end{array}\right.
   \end {equation} in any open set $\Omega\subseteq\mathbb{R}^{3}$ for a.e $ t\in(0,T)$ .

   By the elliptic regularity theory , \begin{equation}
   \|f\|^{\frac{2m}{m+2}}_{L^{\frac{2m}{m+2}}}  \leq\|v\cdot\nabla v\|^{\frac{2m}{m+2}}_{L^{\frac{2m}{m+2}}}+ \|f\|^{\frac{2m}{m+2}}_{H^{-\frac{3}{2}}_{0}}.\end{equation}
So we get
$\nabla p\in L^{\frac{2m}{m+2}}(Q_{T})$ .

   \begin{lem}Let $ u $ be a solution of (1),(2) such that $ u(\cdot,t)$ is analytic in a bounded open set $ Q=\Omega\times(0,T)$,If there exist a nonempty open set $ \Omega_{1}$ in $ \Omega $ and a constant $ t_{1}\in(0,T) $ such that $ u(x,t_{1})=0, x\in \Omega_{1} $,then $ u\equiv 0 $ in $ Q $.\end{lem}

 \textbf{proof}:\quad Since $ u(x,t)$ is analytic in   $ x $  and  $ t $ in $Q$ .By assumption $ u(x,t_{1})=0,$ for $ x\in \Omega_{1}$,hence $ u(x,t_{1})=0$ in $\Omega $ .So $\omega(x,t_{1})=0 $  in  $\Omega $. Since $\omega $ satisfies $$\partial_{t}\omega-\triangle \omega=\ast d (u\nabla u)=div (u\wedge \omega).$$ We have $\partial_{t}\omega(x,t_{1})=0 $ and so $\ast du_{t}(x,t_{1})=0 $ .Since $ u_{t}\in H^{1}_{0}$ and $ div  u_{t} =0 $ we deduce  $ u_{t}(x,t_{1})=0 $.Applying the same argument ,we have $ \frac{\partial}{\partial t}^{k}u(x,t_{1})=0 $ for $ k=0,1,2,...$, then the theorem is proved .

Now we give the new $\varepsilon-$ regularity criterion of Naiver-Stokes equations.
\begin {prop}
For any $\varepsilon>0,$ assume $(u,p)$ satisfies all the conditions of Theorem 1.5 near $(x,t).$ And if
\begin {equation}
r^{\frac{6m}{m+2}-5} \int\int_{Q_{r}(x,t)}|\nabla^{2}u|^{\frac{2m}{m+2}}dxdt\leq r^{\varepsilon},
\end {equation}
for $ r $ is small enough, then (x,t) is a regular point .
\end {prop}
\textbf{proof}: Without loss of generality, we may assume that $(x,t)=(0,0),$ and that $(u,p)$ is defined on a neighborhood $ Q_{\frac{1}{2}}\subseteq D $ of $(0,0)$. First, we note (u,p) satisfies conditions (1.4),(1.5) and is a suitable weak solution to the Navier-Stokes equations in $ D $.

From Lemma 2.2, we obtain
\begin {equation}
\int_{Q_{\frac{1}{4}}}(|\partial_{t}u|^{\frac{2m}{m+2}}+|\nabla^{2}u|^{\frac{2m}{m+2}}+|\nabla p|^{\frac{2m}{m+2}})dz\leq c_{0}
\end {equation}
with an absolute constant $c_{0}.$

Assume that the statement of Proposition 2.4 is false , that is ,$(0,0)$ is a singular point.Then ,as it was shown in [2],there exists a sequence of positive numbers $R_{k}$ such that $R_{k}\rightarrow 0$ as $k\rightarrow\infty$ and
\begin {equation}
A(R_{k})=\frac{1}{R^{2}_{k}}\int_{Q_{R_{k}}}|u(x,t)|^{3}dz>\epsilon_{\ast}
\end {equation}
for all $k\in\mathbb{N}.$ Here $\epsilon_{\ast}$ is an absolute positive constant.
  We extend functions $ u $ and $ p $ to the whole space $\mathbb{R}^{3+1}$ by zero. Extended functions will still be denoted by  $ u $ and $ p ,$ respectively. Now , we let
  $u^{R_{k}}(x,t)=R_{k}u(R_{k}x,R^{2}_{k}t),$
  $p^{R_{k}}(x,t)=R^{2}_{k}p(R_{k}x,R^{2}_{k}t).$

  Obviously,
  \begin {equation}
\int^{1}_{0}(\int_{B_{1}}|u^{R_{k}}(x,t)|^{q}dx)^{\frac{s}{q}}dt
=\frac{1}{R^{\frac{3s}{q}+2-s}}\int^{R^{2}}_{0}(\int_{B_{R}}|u|^{q}dx)^{\frac{s}{q}}dt<C
\end {equation}
\begin {equation}
\int^{1}_{0}\int_{B_{1}}|p^{R_{k}}(x,t)|^{\frac{3}{2}}dxdt=\frac{1}{R^{\frac{7}{2}}}\int^{R^{2}}_{0}\int_{B_{R}}|p|^{\frac{3}{2}}dxdt< C
\end {equation}
for $(q,s)\in\{(2,5),(5,2)\}.$

To extract more information about boundedness of various norms of functions $u^{R_{k}}$ and $p^{R_{k}},$ let us fix a cut-off function $\phi\in C^{\infty}_{0}(\mathbb{R}^{3+1})$ and introduce the function $\phi^{R_{k}}$ in the following way
$$\phi(x,t)=R_{k}\phi^{R_{k}}(R_{k}x,R^{2}_{k}t),x\in\mathbb{R}^{3},t\in\mathbb{R}.$$
We choose $R_{k}$ so small to ensure
$$ supp\phi\subset\{{(x,t)|R^{2}_{k}t\in(-(\frac{1}{4})^{2},(\frac{1}{4})^{2}),R_{k}x\in B(\frac{1}{4}))}\},$$
further,we have
$$ supp\phi^{R_{k}}\subset B(\frac{1}{4}))\times(-(\frac{1}{4})^{2},(\frac{1}{4})^{2}).$$
Then,since the pair$(u,p)$ is a suitable weak solution, we have
$$ 2\int_{Q_{\frac{1}{2}}}\phi^{R_{k}}|\nabla u|^{2}dz\leq\int_{Q_{\frac{1}{2}}}\{|u|^{2}(\triangle\phi^{R_{k}}+\partial_{t}\phi^{R_{k}})+u\cdot\nabla\phi^{R_{k}}(|u|^{2}+2p)\}dz $$
and after changing variable we arrived at the inequality
$$ 2\int_{\mathbb{R}\times\mathbb{R}^{3}}\phi|\nabla u^{R_{k}}|^{2}dz\leq\int_{\mathbb{R}\times\mathbb{R}^{3}}\{|u^{R_{k}}|^{2}(\triangle\phi+\partial_{t}\phi)+u^{R_{k}}
\cdot\nabla\phi(|u^{R_{k}}|^{2}+2p^{R_{k}})\}dz.$$
So, from (2.10),(2.11) and the last two inequalities, we deduce the bound
\begin {equation}
\int_{Q}(|p^{R_{k}}(x,t)|^{\frac{3}{2}}+|\nabla u^{R_{k}}|^{2})dz\leq c_{1}(Q)
\end {equation}
for any domain $Q\Subset\mathbb{R}^{3+1}$ with $c_{1}$ independent of $R_{k}.$ Then we apply known arguments and Lemma 2.2, we find
\begin {equation}
\int_{Q}(|\partial_{t}u^{R_{k}}|^{\frac{2m}{m+2}}+|\nabla^{2}u^{R_{k}}|^{\frac{2m}{m+2}}+|\nabla p^{R_{k}}|^{\frac{2m}{m+2}})dz\leq c_{2}(Q)
\end {equation}
Let us show that
\begin {equation}
u^{R_{k}}\rightarrow v
\end {equation} in $L^{3}(Q_{r})$ for any $0<r<\infty. $
Indeed, by (2.7),
$$\int_{Q_{1}}|\nabla^{2}u^{R_{k}}|^{\frac{2m}{m+2}}dz\leq R^{\varepsilon}_{k}$$
 (2.14) can be easily derived from the interpolation inequality

$$\|u^{R_{k+1}}-u^{R_{k}}\|_{L^{3}(Q_{r})}$$
\begin {equation}
\leq\|u^{R_{k+1}}-u^{R_{k}}\|^{1-\theta}_{L^{\frac{6m}{6-m}}((-r^{2},r^{2}),L^{{\frac{2m}{m+2}}}(B_{r}))}
\|u^{R_{k+1}}-u^{R_{k}}\|^{\theta}_{L^{{\frac{2m}{m+2}}}((-r^{2},r^{2}),L^{\frac{6m}{6-m}}(B_{r}))}.
 \end {equation}

Now,we combine all information about limit $(v,q)$, conclude that:
\begin {equation}
\int_{Q_{\frac{1}{2}}}(|\nabla v|^{2}+|\nabla^{2}v|^{\frac{2m}{m+2}}+|\partial_{t}v|^{\frac{2m}{m+2}}+|\nabla q|^{\frac{2m}{m+2}})dz\leq c_{2}(Q)
\end {equation}
for any $Q\Subset\mathbb{R}^{3+1};$
 \begin {equation}
 \int_{\mathbb{R}\cap Q}(\int_{\mathbb{R}^{3}\cap Q}|v|^{q}dx)^{\frac{s}{q}}ds <c_{2}(Q)
  \end {equation}
  where $ C $ is an absolute constant, and $(q,s)\in\{(2,5),(5,2)\}.$

  Meanwhile, $v$ and $q $ satisfy the Navier-Stokes equations a.e. in $ \mathbb{R}^{3+1},$ that is,
  \begin {equation}
  2\int_{\mathbb{R}}\int_{\mathbb{R}^{3}}\phi|\nabla v|^{2}dz=\int_{\mathbb{R}}\int_{\mathbb{R}^{3}}\{|v|^{2}(\triangle\phi+\partial_{t}\phi)+v\cdot\nabla\phi(|v|^{2}
  +2q)\}dz
  \end {equation}
  for all functions $\phi\in C^{\infty}_{0}(\mathbb{R}^{3+1}).$ It is easy to show that , according to (2.16)-(2.18),the pair
  (v,q) is a suitable weak solution to the Navier-Stokes equations in $\Omega\times[a,b]$ for any bounded domain $\Omega\Subset\mathbb{R}^{3}$ and for any $-\infty<a<b<\infty.$ Moreover, according to (2.9)and (2.14),we find
  \begin {equation}
  \int_{Q\frac{1}{2}}|v|^{3}dz>\epsilon_{\ast}.
  \end {equation}
    Let us proceed the proof of Proposition 2.4. We are going to show there exist some positive numbers $R_{0}$ and $T_{0}$ such that , for any $ k=0,1,\cdots,$ the function $\nabla^{k}v$ is H$\ddot{o}$lder continuous and bounded on the set
    $$(\mathbb{R}^{3}\setminus\bar{B}(\frac{R_{0}}{2}))\times (-T_{0},T_{0}).$$
    To this end , let us fix an arbitrary number $1<T_{0}<2$ and note that
    $$\int^{2T_{0}}_{-2T_{0}}\int_{\mathbb{R}^{3}}(|v|^{3}+|q|^{\frac{3}{2}})dz<\infty.$$
    This means that there exists a number $R_{0}(\varepsilon_{0},T_{0})>2$ such that
    \begin {equation}
    \int^{2T_{0}}_{-2T_{0}}\int_{\mathbb{R}^{3}\setminus\bar{B}(\frac{R_{0}}{2}))}(|v|^{3}+|q|^{\frac{3}{2}})dz<\varepsilon_{0}.
    \end {equation}
    Now,assume that $z_{1}=(x_{1},t_{1})\in(\mathbb{R}^{3}\setminus\bar{B}(R_{0}))\times (-\frac{T_{0}}{2},\frac{T_{0}}{2}).$Then,
    $$Q(z_{1},\frac{1}{2})\equiv B(x_{1},\frac{1}{2})\times(t_{1}-\frac{1}{4},t_{1}+\frac{1}{4})\subset\mathbb{R}^{3}
    \setminus\bar{B}(\frac{R_{0}}{2}))\times(-2T_{0},2T_{0}).$$
    So,by(1.19),
    \begin {equation}
     \int^{t_{1}+\frac{1}{4}}_{t_{1}-\frac{1}{4}}\int_{B(x_{1},\frac{1}{2})}(|v|^{3}+|q|^{\frac{3}{2}})dz<\varepsilon_{0}
     \end {equation}
     for any $z_{1}\in(\mathbb{R}^{3}\setminus\bar{B}(R_{0}))\times (-\frac{T_{0}}{2},\frac{T_{0}}{2}), $ where $1<T_{0}<2$
     and $R_{0}>2.$ Then, it follows from (2.21) and Lemma 1.2, for any $k=0,1,\cdots,$
     \begin {equation}
     \max_{z\in Q(z_{1},\frac{1}{4})}|\nabla^{k}v(z)|\leq c_{0k}<\infty
     \end {equation}
     and $\nabla^{k}v(z)$ is H$\ddot{o}$lder continuous on $(\mathbb{R}^{3}\setminus\bar{B}(R_{0}))\times (-\frac{T_{0}}{2},\frac{T_{0}}{2}).$
     Now,let us introduce the vorticity $\omega$ of $ v $, i.e. $\omega=\nabla\wedge v.$ The function $\omega$ meets the equation
    $$\partial_{t}\omega-\triangle \omega=\ast d (v\cdot\nabla v)=\textsf{div}(v\wedge \omega)$$ in $(\mathbb{R}^{3}\setminus\overline{B}(2R_{0}))\times (-\frac{T_{0}}{4},\frac{T_{0}}{4})$.
    Recalling (2.22), we see that, in the set $(\mathbb{R}^{3}\setminus\overline{B}(2R_{0}))\times (-\frac{T_{0}}{4},\frac{T_{0}}{4}),$ the function $\omega$ satisfies the following relations : \begin{equation}|\partial_{t}\omega-\triangle \omega|\leq M(|\omega|+|\nabla\omega|)        \end{equation} for some constant $ M>0$ and   \begin{equation}|\omega|\leq c_{00}+c_{01}<\infty        .  \end{equation}

    Let us show that
    \begin {equation}
    \omega(z)=0
    \end {equation} for  $a.e.z\in(\mathbb{R}^{3}\setminus\overline{B}(2R_{0}))\times (-\frac{T_{0}}{4},\frac{T_{0}}{4})$
    Indeed,$$(\int_{Q\frac{1}{2}}|v|^{3}dxdt)^{\frac{1}{3}}\leq
    (\int_{Q\frac{1}{2}}|v- u^{R_{k}}|^{3}dxdt)^{\frac{1}{3}}+(\int_{Q\frac{1}{2}}|u^{R_{k}}|^{3}dxdt)^{\frac{1}{3}}$$
    $$\leq(\int_{Q\frac{1}{2}}|v- u^{R_{k}}|^{3}dxdt)^{\frac{1}{3}}
    +\|u^{R_{k}}\|_{L^{\frac{6m}{6-m}}((-\frac{1}{4},\frac{1}{4}),L^{{\frac{2m}{m+2}}}(B_{\frac{1}{2}}))}
    \|u^{R_{k}}\|_{L^{\frac{2m}{m+2}}((-\frac{1}{4},\frac{1}{4}),L^{{\frac{6m}{6-m}}}(B_{\frac{1}{2}}))}$$
    $$\leq(\int_{Q\frac{1}{2}}|v- u^{R_{k}}|^{3}dxdt)^{\frac{1}{3}}
    +\|u^{R_{k}}\|_{L^{\frac{6m}{6-m}}((-\frac{1}{4},\frac{1}{4}),L^{{\frac{2m}{m+2}}}(B_{\frac{1}{2}}))}
    \|\nabla^{2}u^{R_{k}}\|_{L^{\frac{2m}{m+2}}(Q_{\frac{1}{2}})}$$
     $$\leq(\int_{Q\frac{1}{2}}|v- u^{R_{k}}|^{3}dxdt)^{\frac{1}{3}}
    +\|u^{R_{k}}\|_{L^{\frac{6m}{6-m}}((-\frac{1}{4},\frac{1}{4}),L^{{\frac{2m}{m+2}}}(B_{\frac{1}{2}}))}
    \|\nabla^{2}u\|_{L^{\frac{2m}{m+2}}(Q_{\frac{^{R_{k}}}{2}})}\cdot R^{\frac{m-10}{2m}}_{k}.$$
    By (2.7) and (2.14), we can show that
    \begin {equation}
    \int_{Q(\frac{1}{2},z_{\ast})}|v|^{3}dxdt=0
    \end {equation}
    for any $z_{\ast}\in\mathbb{R}^{3+1}.$ So (2.25) is proved.
    Relations (2.23)-(2.25) allow us to apply the backward uniqueness theorem , and conclude that
    \begin {equation}
    \omega(z)=0
    \end {equation} for any $z\in(\mathbb{R}^{3}\setminus\overline{B}(2R_{0}))\times (-\frac{T_{0}}{4},\frac{T_{0}}{4}).$
    If we show that
    \begin {equation}
    \omega(z)=0
    \end {equation}for any $x\in\mathbb{R}^{3},$ a.a.$t\in(-\frac{T_{0}}{4},\frac{T_{0}}{4}),$ then we are done. Indeed, if (2.28) is valid, the function $ v(\cdot,t)$ is harmonic and has the finite $L_{\frac{3}{2}}$-norm. It turn out that  this fact leads to the identity $v(\cdot,t)=0$ for a.a.$t\in(-\frac{T_{0}}{4},\frac{T_{0}}{4}).$ This contradicts with (2.19).

    Now our goal is to show that (2.27) implies (2.28).

      We know that $(v,q)$ meets the equations :  \begin{equation}
      \left\{\begin{array}{l}
      \partial_{t}v+v\cdot\nabla v+\nabla q=0,\\
      \nabla\cdot v=0,\\
      \triangle v=0, \\
      \nabla\wedge v=0
           \end{array}\right.
           \end{equation} in the set $(\mathbb{R}^{3}\setminus\overline{B}(2R_{0}))\times (-\frac{T_{0}}{4},\frac{T_{0}}{4}]$. From (2.29), we deduce the following bound  \begin{equation} \max_{Q_{0}}(|\nabla^{k}v|+|\nabla^{k}\partial_{t}v|+|\nabla^{k}q|\leq c^{1}_{0k}<\infty    \end{equation} for all\quad $k=0,1,2,...$ ,here $Q_{0}=(\mathbb{R}^{3}\setminus\overline{B}(4R_{0}))\times (-\frac{T_{0}}{8},\frac{T_{0}}{8})$.

      To prove (2.28), according to (2.27), we fix a smooth cut-off function $\varphi\in C^{\infty}_{0}(\mathbb{R}^{3})$ subjected to the conditions :$\varphi(x)=1$ if $ x\in B(8R_{0})$ and $\varphi(x)=0 $ if $ x\in\mathbb{R}^{3}\setminus B(12R_{0})$ . Let $ w=\varphi v,r=\varphi q $ ,so $(w,r)$ satisfies
       \begin{equation}
       \left\{\begin{array}{l}
       \partial_{t}w-\triangle w+w\cdot\nabla w+\nabla r=g,\\
        \nabla\cdot w=v\nabla\varphi \end{array}\right.  \end{equation} in $ Q_{\ast}=B(16R_{0})\times (-\frac{T_{0}}{8},\frac{T_{0}}{8})$and \begin{equation} w=0 \quad \texttt{on}\quad\partial B(16R_{0})\times (-\frac{T_{0}}{8},\frac{T_{0}}{8})                  \end{equation} where $ g=(\varphi^{2}-\varphi)v\cdot\nabla v+vv\cdot\nabla\varphi^{2}+q\nabla\varphi-2\nabla v\nabla\varphi-v\triangle\varphi$

It is clear that $w$ is not incompressible. So we introduce the functions $(\widetilde{w},\widetilde{r})$ satisfies :$$-\triangle\widetilde{w}+\nabla\widetilde{r}=0, \nabla\cdot\widetilde{w}= u\nabla\varphi $$ in $ Q_{\ast}$ with $ \widetilde{w}=0 $
on
 $\partial B(16R_{0})\times (-\frac{T_{0}}{8},\frac{T_{0}}{8}).$

 Setting $ U= w-\widetilde{ w} $ and $ P=r-\widetilde{r}$ satisfies  \begin{equation}
       \left\{ \begin{array}{l}
        \partial_{t}U-\triangle U+U\cdot\nabla U+\nabla P=G-div(U\otimes\widetilde{w}+ \widetilde{w}\otimes U),\\
         \nabla\cdot U=0\end{array}\right. \end{equation} in $Q_{\ast}$,  and  \begin{equation} U=0  \quad on \quad  \partial B(16R_{0})\times (-\frac{T_{0}}{8},\frac{T_{0}}{8})    \end{equation}where $ G=-\textsf{div}\widetilde{w}\otimes\widetilde{w}+g-\partial_{t}\widetilde{w} $. By (2.23), (2.16)-(2.18) and the elliptic regularity theory, it lead to the following facts about the smoothness of functions $ U $ and $ P $:
        $$ U \in L^{\infty}((-\frac{T_{0}}{8},\frac{T_{0}}{8});L^{2}(B(16R_{0})))\cap L^{2}((-\frac{T_{0}}{8},\frac{T_{0}}{8});W^{1,2}(B(16R_{0}))),$$
        $$\int_{Q_{\ast}}(|\nabla U|^{2}+|\nabla^{2}U|^{\frac{2m}{m+2}}+|\partial_{t}U|^{\frac{2m}{m+2}}+|\nabla P|^{\frac{2m}{m+2}})dz\leq c_{3}(Q_{\ast}).$$ Furthermore,we can obtain that $(U,P)$ is a suitable weak solution, so the associated space-time singular set S satisfies $\mathfrak{H}^{1}(S)=0.$

         So we can choose $ t_{0}\in(-\frac{T_{0}}{8},\frac{T_{0}}{8})$ so that  \begin{equation}\|\nabla U(\cdot,t_{0})\|_{2, B(16R_{0})}<\infty .\end{equation}
         Then by the short time unique solvability results for the Navier-Stokes equations ,we find a number $\delta_{0}>0 $ such that $$\partial_{t}U,\nabla P,\nabla^{2}U \in L^{2}(B(16R_{0})\times(t_{0},t_{0}+\delta_{0})),$$ then it is easy to check that $$\sup_{t_{0}-\varepsilon<t<t_{0}+\delta_{0}-\varepsilon}\sup_{x\in B(16R_{0})}|\nabla^{k}U|\leq c^{5}_{0k}<\infty $$ for $k=0,1,...$, and for $0<\varepsilon<\frac{\delta_{0}}{4}$, so it is valid that $$\sup_{t_{0}+\varepsilon<t<t_{0}+\delta_{0}-\varepsilon}\sup_{x\in B(16R_{0})}|\nabla^{k}v|\leq c^{6}_{0k}<\infty .$$

          Hence $ v(\cdot,t)$ is analytic in the $ B(8R_{0})$ for $(t_{0}+\varepsilon,t_{0}+\delta_{0}-\varepsilon)$,and as mentioned above, $ \omega=0 $ for  $ (B(8R_{0})\setminus B(4R_{0}))\times(t_{0}+\varepsilon,t_{0}+\delta_{0}-\varepsilon). $  By Lemma 2.3, we obtain $$\omega=0$$ in $ B(8R_{0})\times(t_{0}+\varepsilon,t_{0}+\delta_{0}-\varepsilon)$. Then proposition 2.4 is proved.

          \begin{center}
\item\section{The Main Theorem }
\end{center}

To complete the proof of Theorem 1.5, we give a covering lemma in the following.

          \begin {lem}
          Let $\Re$ be any family of parabolic cylinders $ Q^{\ast}_{r}(x,t)$ contained in a bounded subset of $\mathbb{R}^{3}\times\mathbb{R}.$ Then there exists a finite or denumerable subfamily $\Re'={Q^{\ast}_{i}=Q^{\ast}_{r_{i}}(x_{i},t_{i})}$ such that
          $$ Q^{\ast}_{i}\cap Q^{\ast}_{j}=\varnothing \quad for \quad i\neq j .$$
          $$\forall Q^{\ast}\in\Re,\quad \exists Q^{\ast}_{r_{i}}(x_{i},t_{i})\in\Re',\quad Q^{\ast}\subset Q^{\ast}_{5r_{i}}(x_{i},t_{i}).$$
          \end {lem}

          Now we prove Theorem 1.5: let (u,p) is a weak solution defined on an open set$ D $ and satisfies all conditions in Theorem 1.5; we need to show the conclusion holds in any bounded open set .It is easy to check that (u,p) is a suitable weak solution on any bound open set $\tilde{D}.$ By Proposition 2.4,
          \begin {equation}
        (x,t)\in S \Rightarrow r^{\frac{6m}{m+2}-5} \int\int_{Q_{r}(x,t)}|\nabla^{2}u|^{\frac{2m}{m+2}}dxdt>r^{\varepsilon},
\end {equation}
for $ r $ is small enough.

          Let $ V $ be a neighborhood of $ S $ in  $\tilde{D}.$ And let $\delta>0,$ for each $(x,t)\in S ,$ we choose $Q(x,t)^{\ast}_{r}$ with $ r<\delta $ such that
          $$ r^{\frac{6m}{m+2}-5-\varepsilon} \int\int_{Q^{\ast}_{r}(x,t)}|\nabla^{2}u|^{m}dxdt>C \quad and \quad Q^{\ast}_{r}(x,t)\subset V.$$
            Applying the covering lemma to the family of cylinders, we obtain a disjoint subfamily $ {Q^{\ast}_{r_{i}}(x_{i},t_{i})} $ such that
            $$ S \subset\cup_{i}Q^{\ast}_{5r_{i}}(x_{i},t_{i})$$
            and $$\sum_{i}r^{-\frac{6m}{m+2}+5+\varepsilon}_{i}\leq C^{-1}\sum_{i}\int\int_{Q^{\ast}_{r}}|\nabla^{2}u|^{\frac{2m}{m+2}}dxdt$$
            $$\leq C^{-1}\int\int_{V}|\nabla^{2}u|^{\frac{2m}{m+2}}dxdt .$$
            Since $\delta $ was arbitrary, we conclude that $ S $ has Lebesgue measure zero , and also that
           $$ \mathfrak{H}^{-\frac{6m}{m+2}+5+\varepsilon}(S)\leq\frac{ 5}{C}\int\int_{V}|\nabla^{2}u|^{\frac{2m}{m+2}}dxdt $$
           for every neighborhood $ V $ of $ S .$ Since $|\nabla^{2}u|^{\frac{2m}{m+2}}$ is  integrable, it follows that
           $\mathfrak{H}^{-\frac{6m}{m+2}+5}(S)=0.$

  \end{document}